\documentclass[preprint,amsmath,amssymb,aip]{revtex4-2}

\begin{document}

\title{Classical Nambu brackets in higher dimensions}

\author{Cristel Chandre}
    \affiliation{CNRS, Aix Marseille Univ, I2M, 13009 Marseille, France}
    \email{cristel.chandre@univ-amu.fr}
    
\author{Atsushi Horikoshi}
    \affiliation{Department of Natural Sciences, Tokyo City University,Tokyo 158-8557, Japan}
    \email{horikosi@tcu.ac.jp}


\begin{abstract}
We consider $n$-linear Nambu brackets in dimension $N$ higher than $n$. Starting from a Hamiltonian system with a Poisson bracket and $K$ Casimir invariants defined in the phase space of dimension $N=K+2M$, where $M$ is the number of effective degrees of freedom, we investigate a necessary and sufficient condition for this system to possess $n$-linear Nambu brackets. 
For the case of $n=3$, by looking for the possible solutions to the fundamental identity, the condition is found to be $N=K+2$, i.e., the system should have effectively one degree of freedom. Locally, it is shown that there is only one fundamental solution, up to a local change of variables, and this solution is the canonical Nambu bracket, generated by Levi-Civita tensors. These results generalize to the case of $n(\ge 4)$-linear Nambu brackets.  
\end{abstract}

\maketitle

\section{Introduction}
\label{sec:I}

The Nambu bracket has been introduced for volume-preserving flows in $N$ dimensions as a generalization of Poisson brackets~\cite{Nambu73} (see also Ref.~\onlinecite{Curtright03}).  
The Nambu bracket is an $n$-linear, fully antisymmetric (skew-symmetric) bracket acting on scalar functions (observables) of the variables ${\bf z}\in {\mathbb R}^N$. 
For $n=N=3$ case, it is written as
$$
\{A,B,C\} = 
\epsilon_{ijk}\frac{\partial A}{\partial z_i} \frac{\partial B}{\partial z_j} \frac{\partial C}{\partial z_k},
$$
where $\epsilon$ is the Levi-Civita tensor, $A, B, C$ are 
scalar functions of the variables ${\bf z}\in {\mathbb R}^3$,
and implicit summation (from 1 to $3$) over repeated indices is assumed.
The dynamics of an observable $F({\bf z})$ is generated by two observables $G$ and $H$ from this bracket as
$$
\frac{{\rm d} F}{{\rm d} t} = \{F,G,H\}. 
$$
As a consequence of the antisymmetry, the two observables $G$ and $H$ generating the dynamics are conserved quantities. 
From the identity,
$$
\frac{{\rm d}}{{\rm d}t} \{A, B, C\}=\left\{ \frac{{\rm d} A}{{\rm d} t}, B, C\right\} + \left\{ A, \frac{{\rm d} B}{{\rm d} t}, C\right\} + \left\{ A, B, \frac{{\rm d} C}{{\rm d} t}\right\},
$$
we deduce what is referred to as the fundamental identity:
\begin{equation}
    \label{eqn:FI}
    \{\{A,B,C\},D,E\} = \{\{A,D,E\},B,C\} + \{A, \{B,D,E\},C\} + \{A, B, \{C,D,E\}\}, 
\end{equation}
for all observables $A$, $B$, $C$, $D$ and $E$. In addition, we assume that Nambu brackets satisfy the Leibniz rule, i.e.,
\begin{equation}
    \label{eqn:Leibniz}
\{AB, C, D\} = A\{B, C, D\} + \{A, C, D\} B. 
\end{equation}
The link with Hamiltonian systems comes from defining a bracket
$$
\{A, B\}_C=\{A, B, C\}.
$$
This bracket is a skew-symmetric bi-linear operator, and from the identity~\eqref{eqn:FI} satisfies the Jacobi identity (e.g., by taking $E=C$):
$$
\{\{A, B\}_C, D\}_C + \{\{D, A\}_C, B\}_C + \{\{B, D\}_C, A\}_C =0,
$$
where we have used the antisymmetry of the bracket. 
Therefore, brackets $\{\cdot, \cdot\}_C$ are Poisson brackets for all observables $C$. This is the main argument behind the idea that Nambu brackets are generalizations of Poisson brackets.  

In this article, we investigate properties of skew-symmetric $n$-linear brackets in dimension $N$ higher than $n$, i.e., in phase space spanned by variables ${\bf z}\in {\mathbb R}^N (N\geq n)$. 
The first question here is: Given a Poisson bracket, is it possible to define a Nambu bracket? The answer to this question is rather obvious: 
For example, for the case $n=3$, if a Poisson bracket can be defined from a Nambu bracket, it implies that this Poisson bracket possesses at least one Casimir invariant, namely an observable $C$ which commutes with all the other observables: $\{A,C\}_C = 0$ for all $A$. Obviously, not all Poisson brackets can be derived from a Nambu bracket. So, under which conditions Poisson brackets can be ``paired'' with or derived from Nambu brackets? In addition, what are the different possible Nambu brackets in various dimensions? 
To answer these questions, we consider the case of $n=3$ (and $N\geq 3$) in Sec.~\ref{sec:II},
and extend the results to the case of $n\ge 4$ (and $N\geq n$) in Sec.~\ref{sec:III}. We provide some examples in Sec.~\ref{sec:IV}.

\section{Condition on the Poisson bracket to have a Nambu bracket}
\label{sec:II}

In this section, we consider tri-linear Nambu brackets in dimension $N$ of the form
\begin{equation}
    \label{eqn:Nambu}
    \{A, B, C\} = \lambda_{ijk}({\bf z}) \frac{\partial A}{\partial z_i} \frac{\partial B}{\partial z_j} \frac{\partial C}{\partial z_k},
\end{equation}
where $\lambda_{ijk}$ is a totally antisymmetric tensor called the Nambu tensor, and implicit summation (from 1 to $N$) over repeated indices is assumed. This bracket satisfies the Leibniz rule~\eqref{eqn:Leibniz}, but there are conditions on the Nambu tensor such that the fundamental identity~\eqref{eqn:FI} is satisfied.
More precisely, from the fundamental identity, we deduce two constraints on the Nambu tensor $\lambda$. The first constraint is obtained by looking at terms proportional to $\partial^2 D /\partial z_i\partial z_n$ in Eq.~\eqref{eqn:FI}:
\begin{equation}
    \label{eqn:cond1}
    \lambda_{nij} \lambda_{mkp} + \lambda_{njk} \lambda_{mip} + \lambda_{nki} \lambda_{mjp} +
    \lambda_{mij} \lambda_{nkp} + \lambda_{mjk} \lambda_{nip} + \lambda_{mki} \lambda_{njp} = 0,
\end{equation}
for all $i$, $j$, $k$, $m$, $n$, $p$. We notice that the same condition is obtained by looking at terms proportional to $\partial^2 E /\partial z_i\partial z_n$ (which is not surprising given the antisymmetry of the Nambu bracket, where $D$ and $E$ can easily be exchanged), and that the terms proportional to $\partial^2 A /\partial z_m\partial z_i$, $\partial^2 B /\partial z_m\partial z_i$ and $\partial^2 C /\partial z_m\partial z_i$ cancel out by antisymmetry of the bracket. 
The second constraint on the Nambu tensor $\lambda$ is obtained by looking at terms proportional to the first derivatives of the observables, and writes
\begin{equation}
    \label{eqn:cond2}
    \lambda_{ijk}\frac{\partial \lambda_{mnp}}{\partial z_i} = \lambda_{inp}\frac{\partial \lambda_{mjk}}{\partial z_i} + \lambda_{ipm}\frac{\partial \lambda_{njk}}{\partial z_i} + \lambda_{imn}\frac{\partial \lambda_{pjk}}{\partial z_i},
\end{equation}
for all indices $j$, $k$, $m$, $n$, $p$. We notice the implicit summation over the repeated index $i$. The two conditions~\eqref{eqn:cond1}-\eqref{eqn:cond2} were first derived in Ref.~\onlinecite{Takhtajan94}, and they ensure that the fundamental identity~\eqref{eqn:FI} is satisfied for the Nambu bracket~\eqref{eqn:Nambu}.

As a side note, for a Poisson bracket of the form
$$
\{A, B\} = J_{ij}({\bf z}) \frac{\partial A}{\partial z_i} \frac{\partial B}{\partial z_j},
$$
there exists a well-known condition on the Poisson matrix $\mathbb J$:
\begin{equation}
    \label{eqn:condP}
    J_{il} \frac{\partial J_{jk}}{\partial z_l} + J_{kl} \frac{\partial J_{ij}}{\partial z_l} + J_{jl} \frac{\partial J_{ki}}{\partial z_l} = 0,
\end{equation}
for all indices $i$, $j$, $k$. This corresponds to the condition~\eqref{eqn:cond2} for the Nambu brackets. 
We notice that there are no conditions involving only the values of $J_{ij}$ corresponding to Eq.~\eqref{eqn:cond1}.
As a consequence of this identity, all uniform antisymmetric matrices are Poisson matrices (in the sense that the corresponding bracket is a Poisson bracket) since they satisfy Eq.~\eqref{eqn:condP}. On the contrary, uniform fully-antisymmetric tensors $\lambda$ do not satisfy the constraints imposed by the fundamental identity since Eq.~\eqref{eqn:cond1} is not necessarily satisfied. 

In fact, the problem of finding Nambu tensors $\lambda$ seems to be a much more constrained problem than the one of finding Poisson matrices: Using very rough estimates, a Nambu tensor $\lambda$ is determined by approximately $N^3$ coefficients (less if antisymmetry is taken into account) whereas there are of the order of $N^6$ constraints of the form~\eqref{eqn:cond1} and $N^5$ constraints of the form~\eqref{eqn:cond2}. Again, the number of constraints might by smaller when antisymmetry is taken into account. As a point of comparison, a Poisson matrix is determined by $N^2$ coefficients, and the number of constraints is of order $N^3$. Condition~\eqref{eqn:cond1} seems to be the bottleneck in finding Nambu tensors. 

The objective here is to find all solutions of Eqs.~\eqref{eqn:cond1}-\eqref{eqn:cond2}, and as a consequence all possible Nambu brackets. 
In order to do that, we perform invertible changes of coordinates generated by ${\bf h}: {\mathbb R}^N \to {\mathbb R}^N$, i.e., ${\bf x}={\bf h}({\bf z})$ such that there exists $\widetilde{\bf h}$ such that ${\bf z}=\widetilde{\bf h}({\bf x})$. From the scalar invariance $\widetilde{A}({\bf x}) = A({\bf z})$, the Nambu bracket~\eqref{eqn:Nambu} becomes
$$
\{\widetilde{A}, \widetilde{B}, \widetilde{C}\} = \widetilde{\lambda}_{ijk}({\bf x}) \frac{\partial \widetilde{A}}{\partial x_i} \frac{\partial \widetilde{B}}{\partial x_j} \frac{\partial \widetilde{C}}{\partial x_k},
$$
where 
\begin{equation}
    \label{eqn:equiv}
    \widetilde{\lambda}_{ijk} = \left(\lambda_{mnp}\circ \widetilde{\bf h}\right) \left(\frac{\partial {h}_i}{\partial z_m} \circ \widetilde{\bf h}\right) \left(\frac{\partial {h}_j}{\partial z_n} \circ \widetilde{\bf h}\right) \left(\frac{\partial {h}_k}{\partial z_p} \circ \widetilde{\bf h}\right),
\end{equation}
which can also be written as $\widetilde{\lambda}_{ijk} = \{{h}_i, {h}_j, {h}_k \} \circ \widetilde{\bf h}$. The bracket in the variables ${\bf x}$ is also a Nambu bracket, in the sense that $\widetilde{\lambda}$ satisfies the corresponding conditions~\eqref{eqn:cond1}-\eqref{eqn:cond2} in the variables ${\bf x}$. This can be checked directly or using the identity
$$
\{\widetilde{A}, \widetilde{B}, \widetilde{C}\}({\bf x}) = \{A, B, C\} \circ \widetilde{\bf h}({\bf x}),
$$
in Eq.~\eqref{eqn:FI}. Equation~\eqref{eqn:equiv} defines an equivalence relation in the space of Nambu tensors. Finding all solutions for the Nambu brackets amounts to counting the number of equivalence classes, and identifying a ``canonical'' element in each class.  

Our starting point is a Hamiltonian system in the variables ${\bf z}\in {\mathbb R}^N$ with a Poisson matrix ${\mathbb J}({\bf z})$. As mentioned in Sec.~\ref{sec:I}, if the corresponding Poisson bracket is derived from a Nambu bracket, it means that it has at least one Casimir invariant. This Casimir invariant, denoted $C({\bf z})$ generates the Poisson bracket, i.e., $\{A, B\}_C = \{A, B, C\}$. We assume that the rank of $\mathbb J$ is equal to $N-K$, i.e., that $\mathbb J$ possesses $K$ Casimir invariants. Notice that $N-K$ needs to be even, otherwise, it implies the existence of other Casimir invariants. We denote $N=K+2M$, where $M$ is a positive integer. We then apply the Lie-Darboux theorem~\cite{Arnold78,Libermann87}, which states that locally there exists a change of coordinates ${\bf x} = {\bf h}({\bf z})$ such that the Poisson matrix becomes locally
$$
\widetilde{\mathbb J}({\bf x})=\left( \begin{array}{ccc}
  {\mathbb O}_M   & {\mathbb I}_M & {\mathbb O}_{M, K}  \\
 -{\mathbb I}_M   & {\mathbb O}_M & {\mathbb O}_{M, K}  \\
  {\mathbb O}_{K, M} & {\mathbb O}_{K, M} & {\mathbb O}_{K}   
\end{array} \right),
$$
where ${\mathbb O}_{M, K}$ is the $M\times K$ matrix with all zeros, ${\mathbb O}_{K}$ the $K\times K$ matrix with all zeros, and ${\mathbb I}_M$ the $M\times M$ identity matrix. The Poisson bracket is generated by one of the variables $x_i$ for $i=2M+1,\ldots,N$, let say, $\widetilde{C}({\bf x}) = x_N$ for simplicity in the notations. From $\{\widetilde{A},\widetilde{B},\widetilde{C}\}=\{\widetilde{A},\widetilde{B}\}_{\widetilde{C}}$, we deduce
$$
\widetilde{\lambda}_{ijN} = \widetilde{J}_{ij}.  
$$
Given the antisymmetry of the Nambu tensor, we have
$$
\widetilde{\lambda}_{ijN} = \widetilde{\lambda}_{Nij} = \widetilde{\lambda}_{jNi} = \widetilde{J}_{ij}.
$$
The condition~\eqref{eqn:cond1} for $n=N$ and $p=N$ becomes
$$
\widetilde{J}_{ij} \widetilde{J}_{mk} + \widetilde{J}_{jk} \widetilde{J}_{mi} + \widetilde{J}_{ki} \widetilde{J}_{mj} = 0, 
$$
for all indices $i$, $j$, $k$, $m$. This identity imposes some strong constraints on the possible matrices $\widetilde{\mathbb J}$. For instance, we consider $i=1$, $j=M+1$ and $m=2$, which leads to 
$$
\widetilde{J}_{2,k}=0,
$$
for all $k$ if $M \geq 2$. This results in a contradiction since $\widetilde{J}_{2, M+2}=1$. As a consequence, a necessary condition for the existence of a Nambu bracket is $M=1$. We can check that if $M=1$, the constraints are satisfied for a constant $\widetilde{\lambda}$ defined as
\begin{subequations}
\label{eqn:levi}
\begin{align}
& \widetilde{\lambda}_{ijk} = \epsilon_{ijk}, \;  \mbox{ if } i,j,k \in (1,2,N),\\
& \widetilde{\lambda}_{ijk} = 0, \; \mbox{ otherwise},
\end{align}
\end{subequations}
where $\epsilon$ is the Levi-Civita three-dimensional tensor with indices $(1,2,N)$ instead of the usual $(1,2,3)$. In other words, tri-linear Nambu brackets in higher dimensions are only possible if the number of Casimir invariants of the Poisson bracket is equal to $N-2$, i.e., $N=K+2$. The corresponding dynamical systems are therefore integrable, since they can be reduced to systems with one effective degree of freedom.

{\em Remark}: We have assumed that the Hamiltonian generating the dynamics is time-independent (as a scalar function of the variables ${\bf z}$). If there is a time-dependence in the Hamiltonian, there is a procedure to autonomize the system by including time as an extra variable, together with its canonical conjugate variable. Since the dimension of phase space is now $N+2$, the number of Casimir invariants to get a Nambu bracket is now equal to $N$ which means that the Poisson bracket for the non-autonomous system should be zero, hence its associated dynamics is trivial. In other words, it is not possible to autonomize a system and still get a Nambu bracket in the extended system. 

The question we address now is, given a Poisson bracket generated by Nambu bracket, are there multiple Nambu tensors leading to the same Poisson bracket? To address this point, we work in the variables ${\bf x}$ where locally the Poisson matrix can be rewritten as
\begin{equation}
    \label{eqn:tildeJ}
\widetilde{\mathbb J}({\bf x})=\left( \begin{array}{ccc}
  0   & 1 & {\bf 0}_{K}^\dagger  \\
 -1   & 0 & {\bf 0}_{K}^\dagger  \\
  {\bf 0}_{K} & {\bf 0}_{K} & {\mathbb O}_{K}   
\end{array} \right),
\end{equation}
where $K=N-2$ and ${\bf 0}_{K}$ is a $K$-dimensional vector with zero coefficients. The Poisson bracket becomes
$$
\{\widetilde{A}, \widetilde{B}\} = \frac{\partial \widetilde{A}}{\partial x_1}\frac{\partial \widetilde{B}}{\partial x_2} - \frac{\partial \widetilde{A}}{\partial x_2}\frac{\partial \widetilde{B}}{\partial x_1}.
$$
The question here is how many Nambu tensors can be found such that there exists $\widetilde{C}$ such that $\{\widetilde{A}, \widetilde{B}, \widetilde{C}\} = \{\widetilde{A}, \widetilde{B}\}$ for all observables $\widetilde{A}$ and $\widetilde{B}$. Since $\widetilde{C}$ is a Casimir invariant for the bracket $\{\cdot,\cdot\}$ (by antisymmetry of the bracket $\{\cdot,\cdot,\cdot\}$), it should be a function of $x_k$ for $k=3,\dots,N$. Up to a change of variables, the corresponding Nambu tensor can be reduced to Eq.~\eqref{eqn:levi}. Therefore, up to a change of variables, there is a unique Nambu bracket for each dimension $N$, and it is related to the Levi-Civita tensor as in Eq.~\eqref{eqn:levi}. 

{\em Remark}: As a result of the Lie-Darboux theorem, we notice that there are locally $\lfloor N/2 \rfloor$ possibilities for the Poisson bracket. These possibilities amount to counting the number of pairs of canonically conjugate variables. Concerning the Nambu brackets, there are locally just one bracket for each dimension $N$. This count illustrates that the conditions~\eqref{eqn:cond1}-\eqref{eqn:cond2} resulting from the fundamental identity~\eqref{eqn:FI} are much more restrictive than the condition~\eqref{eqn:condP} resulting from the Jacobi identity. 

{\em Remark}: If we assume that there exists a tri-linear Nambu bracket $\{\cdot,\cdot,\cdot\}$ for a given Hamiltonian system with two invariants, denoted $C({\bf z})$ and $D({\bf z})$, we can define two different Poisson brackets: One is generated by $C$, and the other one is generated by $D$. We denote these two brackets as $\{\cdot, \cdot\}_C = \{\cdot,\cdot, C\}$ and $\{\cdot, \cdot\}_D=\{\cdot,\cdot, D\}$. By applying the fundamental identity~\eqref{eqn:FI}, we get two conditions
\begin{eqnarray*}
&&  \{\{A, B\}_C, E\}_D = \{\{A, E\}_D, B\}_C + \{A, \{B, E\}_D\}_C,\\ 
&& \{\{A, B\}_D, E\}_C = \{\{A, E\}_C, B\}_D + \{A, \{B, E\}_C\}_D,
\end{eqnarray*}
for all observables $A$, $B$, $E$. We now consider a linear combination of these two brackets, i.e., $\{\cdot,\cdot\}_C+\lambda \{\cdot,\cdot\}_D$, where $\lambda$ is a constant. It is rather straightforward to prove the Jacobi identity for this bracket using the two identities above together with the Jacobi identity for each individual bracket (see Ref.~\onlinecite{Esen16}). As a consequence, the Hamiltonian system is bi-Hamiltonian (having two compatible Poisson structures). As a consequence, only bi- (or multi-) Hamiltonian systems can have tri-linear Nambu brackets.

\section{Generalization to multi-linear Nambu brackets}
\label{sec:III}

In this section, we generalize the above results to the case of multi($n\ge 4$)-linear Nambu brackets (see Refs.~\onlinecite{Gautheron96,Grabowski99,Vaisman99,Carinena08} for more details).
We again consider a Hamiltonian system in the variables ${\bf z}\in {\mathbb R}^N$ with a Poisson matrix ${\mathbb J}({\bf z})$. 
We determine the conditions on the Poisson bracket $\{\cdot,\cdot\}$ such that there exists observables $C_k$ for $k=1,\ldots, \kappa$ such that
$$
\{A, B\} = \{A, B, C_1, \ldots, C_\kappa\}= \lambda_{i j k_1\cdots k_\kappa}({\bf z})~ \frac{\partial A}{\partial z_i}\frac{\partial B}{\partial z_j} \frac{\partial C_1}{\partial z_{k_1}}\cdots \frac{\partial C_\kappa}{\partial z_{k_\kappa}},
$$
where the Nambu bracket $\{\cdot,\cdot,\cdot,\ldots,\cdot\}$ is multi-linear (linear in each of its entries), fully antisymmetric (skew-symmetric) and satisfies the Leibniz rule and the fundamental identity (see Ref.~\onlinecite{Takhtajan94}):
\begin{eqnarray}
&& \{\{A, B, C_1, \ldots, C_\kappa\}, D, E_1, \ldots, E_\kappa\} = \{\{A, D, E_1, \ldots, E_\kappa\}, B, C_1, \ldots, C_\kappa\}\nonumber \\
&& \qquad \qquad \qquad + \{A, \{B, D, E_1, \ldots, E_\kappa\}, C_1, \ldots, C_\kappa\} \nonumber \\
&& \qquad \qquad \qquad + \sum_{k=1}^\kappa \{A, B, C_1, \ldots, C_{k-1},\{C_k, D, E_1, \ldots, E_\kappa\}, C_{k+1}, \ldots, C_\kappa\}, \label{eqn:FIdim}
\end{eqnarray}
for all observables $A$, $B$, $C_k$, $D$, $E_k$ for $k=1,\ldots, \kappa$. 
We notice that this implies that the bracket $\{\cdot,\cdot\}$ has at least $\kappa$ Casimir invariants. We assume that ${\mathbb J}({\bf z})$ possesses $K$ Casimir invariants such that $\kappa\leq K$. We consider the tri-linear Nambu bracket
$$
\{A, B, C_1\} = \{A, B, C_1, \ldots, C_\kappa\}. 
$$
By taking $E_k=C_k$ for $k=2,\ldots, \kappa$, Eq.~\eqref{eqn:FIdim} reduces to the fundamental identity~\eqref{eqn:FI}. If a Poisson bracket is derived from a multi-linear Nambu bracket, it is also derived from a tri-linear Nambu bracket. Therefore the condition that a Hamiltonian system has a multi-linear Nambu bracket is that the number of Casimir invariants of the Poisson bracket satisfies $K=N-2$. 

Next let us see the uniqueness of multi-linear Nambu bracket corresponding to the Poisson bracket in the $N=K+2$ system. We assume Nambu brackets of the form
$$
\{A,B,C_1,\ldots, C_\kappa\}= \lambda_{i j k_1\cdots k_\kappa}({\bf z})~ \frac{\partial A}{\partial z_i}\frac{\partial B}{\partial z_j} \frac{\partial C_1}{\partial z_{k_1}}\cdots \frac{\partial C_\kappa}{\partial z_{k_\kappa}},
$$
where $\kappa\le K$ and the $(\kappa+2)$-dimensional tensor is fully antisymmetric. We perform a change of coordinates ${\bf z}\mapsto {\bf x}$ such that $K$ Casimir invariants of the Poisson bracket are $x_3,\ldots, x_N$ and the Poisson matrix is transformed into the form~\eqref{eqn:tildeJ}. We compute 
$$
\{x_1,x_2,x_3,\ldots,x_{\kappa+2}\} = \widetilde{\lambda}_{1,2,3,\ldots, \kappa+2}~({\bf x}) = \{x_1, x_2\} = +1. 
$$
As a consequence, the Nambu tensor is the Levi-Civita tensor of dimension $\kappa+2$. Up to a change of variables, the canonical Nambu brackets generated by Levi-Civita tensors are the only possible bracket satisfying the fundamental identity.

\section{Examples}
\label{sec:IV}

\subsection{Example for $N=3$}
Before proceeding to the higher dimensional $(N\ge 4)$ cases, we consider a case with $N = 3$ given by a Poisson matrix,
$$
{\mathbb J} = \left( 
\begin{array}{ccc}
    0 & 1 & 0  \\
    -1 & 0 & -2q \\
    0 & 2q & 0
\end{array}
\right),
$$
and a Hamiltonian $H({\bf z})$ in three variables ${\bf z}=(q,p,u)$. 
There is one Casimir invariant $C({\bf z})=u-q^2$, and 
the Lie-Darboux form of the Poisson matrix, 
$$
{\widetilde{\mathbb J}_0} = \left( 
\begin{array}{ccc}
    0 & ~1 & ~~0  \\
    -1 & ~0 & ~~0 \\
    0 & ~0 & ~~0
\end{array}
\right),
$$
can be obtained by considering the change of variables ${\bf z}\mapsto {\bf x}=(q,p,C)$.
Since the relation $N=K+2$ is satisfied, there is a tri-linear Nambu bracket with
the canonical Nambu tensor $\lambda_{ijk} = \epsilon_{ijk}$.
From this Nambu bracket and using the Casimir invariant as a Hamiltonian, we easily deduce another Poisson bracket for the same system (see Ref.~\onlinecite{Horikoshi21} for more details).
For example, for a system of a harmonic oscillator with
the Hamiltonian
$$
H({\bf z})=\frac{p^2}{2m} +\alpha u,
$$
we obtain another Poisson structure with the Poisson matrix
$$
\overline{\mathbb J} = \left( 
\begin{array}{ccc}
    0 & ~~\alpha & ~-p/m  \\
    -\alpha & ~~0 & ~0 \\
    p/m & ~~0 & ~0
\end{array}
\right).
$$
Here $-C({\bf z})$ works as a Hamiltonian, 
and $H({\bf z})$ works as a Casimir invariant.
This Poisson matrix $\overline{\mathbb J}$ also reduces to
the canonical symplectic matrix $\widetilde{\mathbb J}_0$
in the variables ${\bf x}=(q/\alpha,p,H)$.

\subsection{Example for $N=4$}
We consider a Hamiltonian system with the Poisson matrix,
$$
{\mathbb J} = \left( 
\begin{array}{cccc}
    0 & 1 & 0 & 2p  \\
    -1 & 0 & -2q & 0\\
    0 & 2q & 0 & 4qp\\
    -2p & 0 & -4qp & 0
\end{array}
\right),
$$
and a Hamiltonian $H({\bf z})$ in four variables ${\bf z}=(q,p,u,v)$. 
There are two Casimir invariants $C_1({\bf z})=u-q^2$ and $C_2({\bf z})=v-p^2$.
The Lie-Darboux form of the Poisson matrix,
$$
\widetilde{\mathbb J}_0 = \left( 
\begin{array}{cccc}
    0 & 1 & ~0 & ~0  \\
    -1 & 0 & ~0 & ~0 \\
    0 & 0 & ~0 & ~0\\
    0 & 0 & ~0 & ~0
\end{array}
\right),
$$
can be obtained by the change of variables ${\bf z}\mapsto {\bf x}=(q,p,C_1,C_2)$, 
Since the relation $N=K+2$ is satisfied there are some noncanonical tri-linear Nambu brackets with Nambu tensors
\begin{eqnarray*}
&& {\lambda}_{ijk} = \epsilon_{ijk}, \;  ~~~~~\mbox{ if } i,j,k \in (1,2,3),\\
&& {\lambda}_{ijk} = -2p\epsilon_{ijk}, \;  \mbox{ if } i,j,k \in (1,3,4),\\
&& {\lambda}_{ijk} = 0, \; ~~~~~~~~\mbox{ otherwise},
\end{eqnarray*}
with $\epsilon_{123}=1$ and $\epsilon_{134}=1$, or 
\begin{eqnarray*}
&& {\lambda}_{ijk} = \epsilon_{ijk}, \;  ~~~~~\mbox{ if } i,j,k \in (1,2,4),\\
&& {\lambda}_{ijk} = -2q\epsilon_{ijk}, \;  \mbox{ if } i,j,k \in (2,3,4),\\
&& {\lambda}_{ijk} = 0, \; ~~~~~~~~\mbox{ otherwise},
\end{eqnarray*}
with $\epsilon_{124}=1$ and $\epsilon_{234}=1$.
Both Nambu tensors reduce to the canonical one: 
\begin{eqnarray*}
&& \widetilde{\lambda}_{ijk} = \epsilon_{ijk}, \;  ~~~~~\mbox{ if } i,j,k \in (1,2,3),\\
&& \widetilde{\lambda}_{ijk} = 0, \; ~~~~~~~~\mbox{ otherwise},
\end{eqnarray*}
with $\epsilon_{123}=1$ in the variables ${\bf x}=(q,p,u,C_2)$ and
${\bf x}=(q,p,v,C_1)$, respectively.
There also exists a canonical quadri-linear Nambu bracket which writes
$$
\{A, B, C, D\} = \epsilon_{ijkl}\frac{\partial A}{\partial z_i} \frac{\partial B}{\partial z_j}\frac{\partial C}{\partial z_k}\frac{\partial D}{\partial z_l}, 
$$
where $\epsilon$ is the four-dimensional Levi-Civita tensor. 
From this quadri-linear Nambu bracket, the above two tri-linear Nambu brackets can be deduced as $\{A, B, C, C_2\}=\{A, B, C\}_{C_2}$ and 
$\{A, B, C_1, C\}=\{A, B, C\}_{C_1}$, respectively. 

If the Hamiltonian is given, another noncanonical tri-linear Nambu bracket can be obtained as $\{A, H, B, C\}=\{A, B, C\}_{H}$.
For example, for a system of a harmonic oscillator with
the Hamiltonian
$$
H({\bf z})=\frac{v}{2m} +\alpha u,
$$
we obtain another Nambu bracket with the Nambu tensor
\begin{eqnarray*}
&& {\lambda}_{ijk} = 1/(2m)\epsilon_{ijk}, \;  ~\mbox{ if } i,j,k \in (1,2,3),\\
&& {\lambda}_{ijk} = -2\alpha\epsilon_{ijk}, \;  ~~~~~\mbox{ if } i,j,k \in (1,2,4),\\
&& {\lambda}_{ijk} = 0, \; ~~~~~~~~~~~~~~\mbox{ otherwise},
\end{eqnarray*}
with $\epsilon_{123}=1$ and $\epsilon_{124}=1$.
This Nambu tensor also reduces to the canonical one in the variables 
${\bf x}=(q,2mp,u,H)$. 
From these tri-linear Nambu brackets, we deduce two other Poisson structures:
$$
{\mathbb J}_1 = \left( 
\begin{array}{cccc}
    0 & ~~\alpha & ~-p/m & ~2\alpha p  \\
    -\alpha & ~~0 & ~~0 & ~0\\
    p/m & ~~0 & ~~0 & ~0\\
    -2\alpha p & ~~0 & ~~0 & ~0
\end{array}
\right),
$$
for the Hamiltonian $-C_1$ (and $H$ and $C_2$ as Casimir invariants), and
$$
{\mathbb J}_2 = \left( 
\begin{array}{cccc}
    0 & 1/(2m) & 0 & ~~0  \\
    -1/(2m) & 0 & -q/m & ~2\alpha q\\
    0 & q/m & 0 & ~~0\\
    0 & -2\alpha q & 0 & ~~0
\end{array}
\right),
$$
for the Hamiltonian $-C_2$ (and $H$ and $C_1$ as Casimir invariants). The Poisson matrices ${\mathbb J}_1$ and ${\mathbb J}_2$ reduce locally to the canonical symplectic matrix $\widetilde{\mathbb J}_0$
in the variables ${\bf x}=(q/\alpha,p,H,C_2)$ and ${\bf x}=(q,2mp,C_1,H)$, respectively. 

\subsection{Example for $N=6$}
We now consider a system of six variables ${\bf z}=(q_1, p_1, u_1, q_2, p_2, u_2)$ with the Poisson matrix,
$$
{\mathbb J} = \left( 
\begin{array}{cccccc}
  0 & ~1 & 0 & 0 & ~0 & 0  \\
  -1 & ~0 & -2q_1 & 0 & ~0 & 0 \\
  0 & ~2q_1 & 0 & 0 & ~0 & 0 \\
  0 & ~0 & 0 & 0 & ~1 & 0 \\
  0 & ~0 & 0 & -1 & ~0 & -2q_2 \\
  0 & ~0 & 0 & 0 & ~2q_2 & 0
\end{array}
\right),
$$
and a Hamiltonian $H({\bf z})$~\cite{Horikoshi21}.
There are two Casimir invariants $C_1({\bf z})=u_1-q_1^2$ and $C_2({\bf z})=u_2-q_2^2$.
The Lie-Darboux form of the Poisson matrix can be obtained by the change of variables ${\bf z}\mapsto {\bf x}=(q_1, q_2, p_1, p_2, C_1, C_2)$,
$$
\widetilde{\mathbb J}_0 = \left( 
\begin{array}{cccccc}
  0 & 0 & ~1 & ~~0 & ~~0 & ~~0  \\
  0 & 0 & ~0 & ~~1 & ~~0 & ~~0 \\
  -1 & 0 & ~0 & ~~0 & ~~0 & ~~0 \\
  0 & -1 & ~0 & ~~0 & ~~0 & ~~0 \\
  0 & 0 & ~0 & ~~0 & ~~0 & ~~0\\
  0 & 0 & ~0 & ~~0 & ~~0 & ~~0
\end{array}
\right).
$$
Since $N=6$ and $K=2$, the relation $N=K+2$ is not satisfied and 
we conclude that a Nambu bracket for this system does not exist, whether it is a tri- or quadri-linear bracket. 
For example, consider a three-dimensional tensor,
\begin{eqnarray*}
&& {\lambda}_{ijk} = \epsilon_{ijk}, \;  ~\mbox{ if } i,j,k \in (1,2,3),\\
&& {\lambda}_{ijk} = \epsilon_{ijk}, \;  ~\mbox{ if } i,j,k \in (4,5,6),\\
&& {\lambda}_{ijk} = 0, \; ~~~~\mbox{ otherwise},
\end{eqnarray*}
with $\epsilon_{123}=1$ and $\epsilon_{456}=1$.
Using this tensor we can derive the Poisson structure ${\mathbb J}$,
however, this tensor does not satisfy the condition~\eqref{eqn:cond1}.
 
We notice that the tri-linear bracket defined by the above three-dimensional tensor $\lambda_{ijk}$ satisfies the fundamental identity~\eqref{eqn:FI}
if the two degrees of freedom $(q_1,p_1,u_1)$ and $(q_2,p_2,u_2)$ are decoupled in the functions $D$ and $E$ in Eq.~\eqref{eqn:FI} (see Ref.~\onlinecite{Horikoshi21} for more details).
For example, we consider a system of interacting two oscillators 
whose Hamiltonian is given by 
$$
H({\bf z})=\frac{p_1^2}{2m_1} + \frac{p_2^2}{2m_2} +\alpha_1 u_1 +\alpha_2 u_2 + \mu q_1 u_2,
$$
and choose $D=H$ and $E=C_1+C_2$ in the fundamental identity~\eqref{eqn:FI}.
If the interaction is turned off, $\mu=0$, the identity holds, whereas if the interaction is turned on, $\mu\ne0$, the identity no longer holds.
Given that the system reduces to a simplified H\'enon-Heiles system 
$$
H(q_1,q_2,p_1,p_2)=\frac{p_1^2}{2m_1}+\frac{p_2^2}{2m_2}+\alpha_1 q_1^2+\alpha_2 q_2^2+\mu q_1 q_2^2,
$$
its dynamics has two degrees of freedom and is not integrable for $\mu\ne0$~\cite{Heller80}, a further indication that a Nambu bracket for this system does not exist. 

As a side note, as in the case for $N=4$, there is a Poisson structure for this system with $-C_1$ as Hamiltonian (and $H$ and $C_2$ as Casimir invariants).
In order to derive it, it is easier to move into a coordinate system $(q_1,p_1,q_2,p_2,H, C_2)$, and in these variables, it is given by
$$
\widetilde{\mathbb J}_1 = \left( 
\begin{array}{cccccc}
  0 & \alpha_1 & ~0 & ~0 & ~~0 & ~~0  \\
  -\alpha_1 & 0 & ~0 & ~0 & ~~0 & ~~0 \\
  0 & 0 & ~0 & ~\alpha_1 & ~~0 & ~~0 \\
  0 & 0 & -\alpha_1 & ~0 & ~~0 & ~~0 \\
  0 & 0 & ~0 & ~0 & ~~0 & ~~0\\
  0 & 0 & ~0 & ~0 & ~~0 & ~~0
\end{array}
\right),
$$
and it is straightforward to check the Jacobi identity in these variables. In the original variables ${\bf z}$, the Poisson matrix writes
$$
{\mathbb J}_1 = \left( 
\begin{array}{cccccc}
  0 & ~\alpha_1 & -p_1/m_1 & ~~0 & ~~0 & ~~0  \\
  -\alpha_1 & ~0 & \mu u_2 & ~~0 & ~~0 & ~~0 \\
  p_1/m_1 & ~-\mu u_2 & 0 & ~p_2/m_2 & -2q_2(\alpha_2+\mu q_1) & 2q_2p_2/m_2 \\
  0 & ~0 & -p_2/m_2 & ~~0 & ~~\alpha_1 & ~~0 \\
  0 & ~0 & 2q_2(\alpha_2+\mu q_1) & ~-\alpha_1 & ~~0 & -2\alpha_1 q_2\\
  0 & ~0 & -2q_2p_2/m_2 & ~~0 & ~2\alpha_1 q_2 & ~~0
\end{array}
\right).
$$
It can also be checked that ${\mathbb J}+\lambda {\mathbb J}_1$, where $\lambda $ is constant, is also a Poisson bracket. Consequently, the system with $N=6$ possesses two compatible Poisson structures, and is therefore bi-Hamiltonian. We conclude that any system which possess a Nambu bracket is bi- (or multi-) Hamiltonian, but the converse is not true: There exists bi-Hamiltonian systems which do not possess a Nambu bracket. In finite dimension, Nambu systems form a subset of bi-Hamiltonian systems.

\section*{Conclusions}

For a given Poisson bracket $\{\cdot, \cdot \}$ with a Casimir invariant $C$, is there a tri-linear Nambu bracket $\{\cdot, \cdot,\cdot \}$ such that $\{A,B\}=\{A,B,C\}$ for all observables $A$ and $B$? Here we replied that this is only feasible if the number of Casimir invariants $K$ and the dimension of phase space $N$ satisfy the relation $N=K+2$. 
That is, there is a Nambu bracket corresponding to the Poisson bracket only in a Hamiltonian system with one effective degree of freedom.
In this case, one of our conclusions is that locally, up to a change of variables, all Nambu brackets are canonical, in the sense that it is trivially linked to the Levi-Civita tensor. The constraints generated by the fundamental identity are so strong that basically there is only one fundamental solution for the Nambu bracket, up to a local change of coordinates. This conclusion is generalized to multi-linear Nambu brackets. 

\section*{Data Availability}
Data sharing is not applicable to this article as no new data were created or analyzed in this study.


%

\end{document}